\documentclass[11pt]{amsart}
\usepackage{epsfig}
\usepackage{graphicx}
\usepackage{amscd}
\usepackage{amsmath}
\usepackage{amsxtra}
\usepackage{amsfonts}
\usepackage{amssymb}

\newtheorem{theorem}{Theorem}[section]
\newtheorem{corollary}[theorem]{Corollary}
\newtheorem{lemma}[theorem]{Lemma}
\newtheorem{proposition}[theorem]{Proposition}

\theoremstyle{definition}
\newtheorem{definition}[theorem]{Definition}
\newtheorem{remark}[theorem]{Remark}

\newtheorem{example}[theorem]{Example}
\theoremstyle{remark}

\renewcommand{\theclaim}{\textup{\theclaim}}

\numberwithin{equation}{section}

\def\openone

{\mathchoice

{\hbox{\upshape \small1\kern-3.3pt\normalsize1}}

{\hbox{\upshape \small1\kern-3.3pt\normalsize1}}

{\hbox{\upshape \tiny1\kern-2.3pt\SMALL1}}

{\hbox{\upshape \Tiny1\kern-2pt\tiny1}}}

\makeatletter

\newbox\ipbox

\newcommand{\ip}[2]{\left\langle #1\, , \,#2\right\rangle}
\newcommand{\diracb}[1]{\left\langle #1\mathrel{\mathchoice

{\setbox\ipbox=\hbox{$\displaystyle \left\langle\mathstrut
#1\right.$}

\vrule height\ht\ipbox width0.25pt depth\dp\ipbox}

{\setbox\ipbox=\hbox{$\textstyle \left\langle\mathstrut
#1\right.$}

\vrule height\ht\ipbox width0.25pt depth\dp\ipbox}

{\setbox\ipbox=\hbox{$\scriptstyle \left\langle\mathstrut
#1\right.$}

\vrule height\ht\ipbox width0.25pt depth\dp\ipbox}

{\setbox\ipbox=\hbox{$\scriptscriptstyle \left\langle\mathstrut
#1\right.$}

\vrule height\ht\ipbox width0.25pt depth\dp\ipbox}

}\right. }

\newcommand{\dirack}[1]{\left. \mathrel{\mathchoice

{\setbox\ipbox=\hbox{$\displaystyle \left.\mathstrut
#1\right\rangle$}

\vrule height\ht\ipbox width0.25pt depth\dp\ipbox}

{\setbox\ipbox=\hbox{$\textstyle \left.\mathstrut
#1\right\rangle$}

\vrule height\ht\ipbox width0.25pt depth\dp\ipbox}

{\setbox\ipbox=\hbox{$\scriptstyle \left.\mathstrut
#1\right\rangle$}

\vrule height\ht\ipbox width0.25pt depth\dp\ipbox}

{\setbox\ipbox=\hbox{$\scriptscriptstyle \left.\mathstrut
#1\right\rangle$}

\vrule height\ht\ipbox width0.25pt depth\dp\ipbox}

} #1\right\rangle}

\newcommand{\cj}[1]{\overline{#1}}

\newcommand{\bz}{\mathbb{Z}}

\newcommand{\br}{\mathbb{R}}

\newcommand{\bn}{\mathbb{N}}

\def\blfootnote{\xdef\@thefnmark{}\@footnotetext}

\renewcommand{\mod}{\operatorname{mod}}

\newcommand{\dist}{\operatorname*{dist}}
\newcommand{\clspan}{\overline{\operatorname*{span}}}
\input xy
\xyoption{all}
\usepackage{amssymb}

\newcommand{\ket}[1]{\lvert #1\rangle}
\newcommand{\bra}[1]{\langle #1\rvert}

\def\R{\mathcal{R}}

\def\H{\mathcal{H}}

\def\-{^{-1}}

\def\K{\mathcal{K}}

\def\Rs{{R^T}}
\def\Rsi{(R^T)^{-1}}

\begin{document}

\title[Spectral measures and Cuntz algebras]{Spectral measures and Cuntz algebras}
\author{Dorin Ervin Dutkay}
\blfootnote{}
\address{[Dorin Ervin Dutkay] University of Central Florida\\
	Department of Mathematics\\
	4000 Central Florida Blvd.\\
	P.O. Box 161364\\
	Orlando, FL 32816-1364\\
U.S.A.\\} \email{ddutkay@mail.ucf.edu}

\author{Palle E.T. Jorgensen}
\address{[Palle E.T. Jorgensen]University of Iowa\\
Department of Mathematics\\
14 MacLean Hall\\
Iowa City, IA 52242-1419\\}\email{jorgen@math.uiowa.edu}

\thanks{With partial support by the National Science Foundation} 

\begin{abstract}
  We consider a family of measures  $\mu$ supported in $\br^d$ and generated in the sense of Hutchinson by a finite family of affine transformations. It is known that interesting sub-families of these measures allow for an orthogonal basis in $L^2(\mu)$ consisting of complex exponentials, i.e., a Fourier basis corresponding to a discrete subset $\Gamma$ in $\br^d$. Here we offer two computational devices for understanding the interplay between the possibilities for such sets $\Gamma$ (spectrum) and the measures $\mu$ themselves. Our computations combine the following three tools: duality, discrete harmonic analysis, and dynamical systems based on representations of the Cuntz $C^*$-algebras $\mathcal O_N$. 
\end{abstract}
\maketitle \tableofcontents
\section{Introduction}\label{intr}

      The idea of selfsimilarity is ubiquitous in pure and applied mathematics and has many incarnations, and even more applications. In this generality the method often goes by the name multi-scale theory, see for example \cite{ZHS09}. We are faced with a fixed system $S$, say finite, of transformations. The mappings from $S$ are applied to the data at hand, and we look for similarity. A single measure, say $\mu$, may be transformed with the mappings in $S$. If $\mu$ is a convex combination of the resulting measures, we say that $\mu$ is a selfsimilar measure. Here we aim for a harmonic analysis of the Hilbert space $L^2(\mu)$ in the event $\mu$ is selfsimilar. However this is meaningful only if further restrictions are placed on the setup. For example, we will assume that $S$ consists of affine and contractive transformations in $\br^d$ for some fixed dimension $d$; and we will assume that the affine transformations arise from scaling with the same matrix for the different mappings in $S$. In this case $\mu$ is the result of an iteration in the small, and its support is a compact subset in $\br^d$. We explore algorithms for generating Fourier bases of complex exponentials in $L^2(\mu)$ by iteration in the large.

      A word about the Cuntz algebras $\mathcal O_N$ from the title of our paper: These are
infinite algebras on a finite number of generators, the following relations
\begin{equation}\label{eqcuntz}
S_i^*S_j=\delta_{ij}1,\quad(i,j\in\{1,\dots,N\}),\quad\sum_{i=1}^NS_iS_i^*=1.
\end{equation}
 They are intrinsically selfsimilar
and therefore ideally serve to encode iterated function systems (IFSs) from
geometric analysis, and their measures $\mu$. At the same time, their
representations offer (in a more subtle way) a new harmonic analysis of the
associated $L^2(\mu)$-Hilbert spaces. Even though the Cuntz algebras initially
entered into the study of $C^*$-algebras \cite{Cun77} and physics, in recent years
these same Cuntz algebras, and their representation, have found increasing
use in both pure and applied problems, wavelets, fractals, signals; see for
example \cite{Jor06}.

       Earlier work along these lines include the papers \cite{DuJo06b, DJ07a, Fu74, HeLa08, Jor82, Ped04, Ta04}; as well as an array of diverse applications  \cite{ACH07, BaVi05, CoMa07, DeHe09, HuLa08, MoZu09, OkSt05, YAL03}. Recent papers by Barnsley et al  \cite{Bar09, BHS05} show that even the restricted family of affine IFSs suffices for such varied applications as vision, music, Gabor transform processing, and Monte Carlo algorithms. The literature is large and we direct the reader to \cite{DJ07a} for additional references. Our present study suggests algorithms for this problem based on representations of a certain scale of $C^*$-algebras (the Cuntz algebras \cite{Cun77}); see also \cite{Jor06}. In addition we shall make use of Hadamard matrices \cite{ALM02, CHK97, Di04, GoRo09, KhSe93}; number theory and tilings \cite{BaLu07, FLS09, LeLa07, MuPo04}; and of tools from symbolic dynamics \cite{Den09, Tho05, DFG04, MiMi09}.

 Our starting point is a system $(R, B, L)$ in $\br^d$ subject to the conditions from Definitions \ref{def1.2}, i.e., a Hadamard triple. With the given expansive matrix $R$ and the vectors from the finite set $B$ we build one affine iterated function system (IFS) \eqref{eq1.1}, and the associated Hutchinson measure $\mu_B$. Our objective is to set up an algorithmic approach for constructing orthogonal families of complex exponentials in $L^2(\mu_B)$, so called Fourier bases.

    With the third part of the triple $(R, B, L)$, the finite set of vectors from $L$ and the transposed matrix $R^T$ we build a second IFS \eqref{eq1.12}, the $L$-system, now with $R^T$ as scaling matrix and the vectors from $L$ as translations. The set $L$, together with the transpose matrix $R^T$ will be used to construct the frequencies associated to the Fourier basis. Due to the Hadamard property which is assumed for the combined system $(R, B, L)$, we then get a natural representation of the Cuntz algebra $\mathcal O_N$ (where $N = \# L = \# B$) with the generating operators $\{S_l : l \in L\}$ acting as isometries in $L^2(\mu_B)$, see \eqref{eq1.9} in Proposition \ref{pr1.3}.

     Our objective is to use this Cuntz-algebra representation in order to recursively construct and exhibit an orthogonal family of complex exponentials in $L^2(\mu_B)$. For this purpose we introduce a family of minimal and finite cycles $C$ for the $L$-system, which are extreme in the sense made precise in Definition \ref{defcyc}. For each of the extreme cycles $C$, with the $L$-Cuntz representation, we then generate an infinite orthogonal family $\Gamma(C)$ of complex exponentials in $L^2(\mu_B)$. With the cycle $C$ fixed, we prove that the corresponding closed subspace $\H(C)$ now reduces the Cuntz-algebra representation, and further that the restricted representation is irreducible. For distinct minimal cycles we get orthogonal subspaces $\H(C)$ and corresponding disjoint representations (see Theorem \ref{th1.8}).

    A key tool in the study of these representations of the Cuntz algebra is a transfer operator $R_{B,L}$ (Definition \ref{deftro}), analogous to a transfer operator used first by David Ruelle in a different context. In Corollary \ref{corcto} we show that the reduction of the $L$-Cuntz representation \eqref{eq1.9} is accounted for by a family of harmonic functions $h_C$ for the transfer operator $R_{B,L}$.  

  Finally we prove that the sum of the functions $h_C$ equals to the constant function 1 if and only if the union of the sets $\Gamma(C)$ forms an orthogonal basis of Fourier exponentials in $L^2(\mu_B)$.

  In section 1, we define the representation of the Cuntz algebra $(S_l)_{l\in L}$ associated to a Hadamard pair and we present some computational features of this representation. We show in Proposition \ref{prrht} and Corollary \ref{corcto} how the canonical endomorphism on $\mathcal B(L^2(\mu_B))$, constructed from the representation:
$$\alpha(T)=\sum_{l\in L}S_lTS_l^*,\quad(T\in\mathcal B(L^2(\mu_B)))$$
is connected to the transfer operator associated to the ``dual'' IFS $(\tau_l)_{l\in L}$. 

 In section 2 we revisit the permutative representations defined in \cite{BrJo99}. They are needed for the decomposition of our representation $(S_l)_{l\in L}$ into irreducible representations. We will see in section 3, that extreme cycles (Definition \ref{defcyc}) will generate irreducible atoms in this decomposition (Theorem \ref{th1.8}). In dimension 1, we know from \cite{DJ06b} that the extreme cycles offer a complete description of the picture. However, in higher dimensions, more complicated decompositions might appear. We analyze these possibilities in section 4, and we illustrate it in Example \ref{ex4.7}.

\begin{definition}\label{def1.3}

We will denote by $e_t$ the exponential function
$$e_t(x)=e^{2\pi it\cdot x},\quad(x,t\in\br^d)$$

Let $\mu$ be a Borel probability measure on $\br^d$. We say that $\mu$ is a {\it spectral} measure if there exists a subset $\Lambda$ of $\br^d$ such that the family $E(\Lambda):=\{e_\lambda : \lambda\in\Lambda\}$ is an orthonormal basis for $L^2(\mu)$. In this case $\Lambda$ is called a {\it spectrum} for the measure $\mu$ and we say that $(\mu,\Lambda)$ is a{\it spectral pair}.
\end{definition}

\begin{definition}\label{def1.1}
We will use the following assumptions throughout the paper.
Let $R$ be a $d\times d$ expansive integer matrix, i.e., all its eigenvalues have absolute value strictly bigger than one, and let $B$ be a finite subset of $\bz^d$, with $0\in B$. We denote by $N$ the cardinality of $B$. Define the maps
\begin{equation}
\tau_b(x)=R^{-1}(x+b),\quad(x\in\br^d,b\in B)
\label{eq1.1}
\end{equation}

We call $(\tau_b)_{b\in B}$ the (affine) iterated function system (IFS) associated to $R$ and $B$. 

By \cite{Hut81}, there exists a unique compact set $X_B$ called {\it the attractor} of the IFS $(\tau_b)_{b\in B}$ such that
\begin{equation}
X_B=\bigcup_{b\in B}\tau_b(X_B).
\label{eq1.2}
\end{equation}
In our case, it can be written explicitly
\begin{equation}
X_B=\left\{\sum_{k=1}^\infty R^{-k}b_k : b_k\in B\mbox{ for all }k\geq1\right\}.
\label{eq1.3}
\end{equation}
We will use also the notation $X(B)$ for $X_B$.

There exists a unique Borel probability measure $\mu_B$ such that 
\begin{equation}
\mu_B(E)=\frac{1}{N}\sum_{b\in B}\mu_B(\tau_b^{-1}(E))\mbox{ for all Borel subsets }E\mbox{ of }\br^d.
\label{eq1.4}
\end{equation}
Equivalently 
\begin{equation}
\int f\,d\mu_B=\frac1N\sum_{b \in B}\int f\circ \tau_b\,d\mu_B\mbox{ for all bounded Borel function }f\mbox{ on }\br^d.
\label{eq1.5}
\end{equation}
The measure $\mu_B$ is called the {\it invariant measure} of the IFS $(\tau_b)_{b\in B}$. It is supported on $X_B$.

We say that the measure $\mu_B$ {\it has no overlap} if 
\begin{equation}
\mu_B(\tau_b(X_B)\cap \tau_{b'}(X_B))=0,\mbox{ for all }b\neq b'\in B.
\label{eq1.6}
\end{equation}

If the measure $\mu_B$ has no overlap, then one can define the map $\R:X_B\rightarrow X_B$
\begin{equation}
\R(x)=Rx-b,\mbox{ if }x\in\tau_b(X_B).
\label{eq1.7}
\end{equation}
The map $\R$ is well defined $\mu_B$-a.e. on $X_B$.
\end{definition}

\begin{definition}\label{def1.2}
Let $R$ be a $d\times d$ integer matrix, and $B$ and $L$ two subsets of $\bz^d$ of the same cardinality $N$, with $0\in B$ and $0\in L$. We say that $(B, L)$ forms a {\it Hadamard pair} if the following matrix is unitary:
\begin{equation}
\frac{1}{\sqrt N}\left(e^{2\pi i R^{-1}b\cdot l}\right)_{b\in B,l\in L}.
\label{eq1.8}
\end{equation}

\end{definition}

Throughout, the Hilbert space considered is $L^2(\mu_B)$ unless otherwise
specified. And we introduce specific families of operators acting there,
starting with the operator system $(S_l)_{l\in L}$ in \eqref{eq1.9} below. This system $(S_l)_{l\in L}$
will be fixed, and it defines a representation of the Cuntz algebra $\mathcal O_N$ with $N=\#L$.

\begin{proposition}\label{pr1.3}
Let $(\tau_b)_{b\in B}$ be as in Definition\ref{def1.1} and assume the invariant measur $\mu_B$ has no overlap. Let $L$ be a subset of $\bz^d$ with the same cardinality $N$ as $B$. Define the operators on $L^2(\mu_B)$:
\begin{equation}
(S_lf)(x)=e^{2\pi i l\cdot x}f(\R x),\quad(x\in X_B,f\in L^2(\mu_B))
\label{eq1.9}
\end{equation} 
\begin{enumerate}
\item
The operator $S_l$ is an isometry for all $l\in L$.
\item
The operators $(S_l)_{l\in L}$ form a representation of the Cuntz algebra $\mathcal O_N$ if and only if $(B,L)$ forms a Hadamard pair. 

\end{enumerate}

\end{proposition}

\begin{proof}
We will need the following lemma:
\begin{lemma}\label{lem1.3}
If the measure $\mu_B$ has no overlap then for all $b\in B$ and all integrable Borel functions $f$:
\begin{equation}
\int_{\tau_b(X_B)} f\,d\mu_B=\frac1N\int f\circ\tau_b\,d\mu_B
\label{eq1.10}
\end{equation}
\end{lemma}
\begin{proof}
From the invariance equation
$$N\int\chi_{\tau_{b}(X_B)}g\,d\mu_B=\sum_{b'\in B}\int\chi_{\tau_{b}(X_B)}\circ\tau_b'\,g\circ\tau_b'\,d\mu_B.$$
Since there is no overlap $\chi_{\tau_{b}(X_B)}\circ\tau_{b'}$ is 1 if $b=b'$, and is 0 if $b\neq b'$, $\mu_B$-almost everywhere (since $\mu_B$ is supported on $X_B$).

This leads to the conclusion.
\end{proof}

We prove (i). Take $f\in L^2(\mu_B)$, $l\in L$. Then, with Lemma \ref{lem1.3}, 
$$\int |e^{2\pi il\cdot x}f(\R x)|^2\,d\mu_B=\sum_{b\in B}\int_{\tau_b(X_B)}|f(Rx-b)|^2\,d\mu_B(x)$$$$=
\sum_{b\in B}\frac1N\int |f(R\tau_b(x)-b)|^2\,d\mu_B(x)=\int |f|^2\,d\mu_B.$$
This shows that $S_l$ is an isometry.

We prove now (ii). For this, we have to compute $S_l^*$. We have, for $f,g\in L^2(\mu_B)$:
$$\ip{S_lf}{g}=\int e^{2\pi il\cdot x}f(\R x)\cj g(x)\,d\mu_B(x)=
\sum_{b\in B}\int_{\tau_b(X)}e^{2\pi i\cdot x}f(Rx-b)\cj g(x)\,d\mu_B(x)\mbox{, and with Lemma \ref{lem1.3},}$$
$$=\frac{1}{N}\sum_{b\in B}\int e^{2\pi i\cdot\tau_b(x)}f(x) \cj g(\tau_b(x))\,d\mu_B(x).$$
This shows that 
\begin{equation}
(S_l^*g)(x)=\frac{1}{N}\sum_{b\in B}e^{-2\pi i l\cdot\tau_b(x)}g(\tau_b(x)),\quad(x\in X_B,g\in L^2(\mu_B))
\label{eq1.11}
\end{equation}

Then
$$S_l^*S_{l'}f(x)=\frac1N\sum_{b\in B}e^{-2\pi il\cdot \tau_b(x)}e^{2\pi i l'\tau_b(x)}f(R\tau_b(x)-b)=
\left(\frac{1}{N}\sum_{b\in B}e^{2\pi i(l'-l)\cdot R^{-1}b}\right)e^{2\pi i(l'-l)\cdot R^{-1}x}f(x).$$

Therefore $S_l^*S_{l'}=\delta_{l,l'}I_{L^2(\mu_B)}$ if and only if the matrix in \eqref{eq1.8} is unitary.

Also 
$$A:=\sum_{l\in L}S_lS_l^*f(x)=\sum_{l\in L}e^{2\pi il \cdot  x}\frac{1}{N}\sum_{b\in B}e^{-2\pi i l\cdot\tau_b(\R x)}f(\tau_b(\R x)).$$
If $x\in\tau_{b_0}(X_B)$, then $\tau_b(\R x)=x+b-b_0$ and we have further
$$A=\sum_{b\in B}f(x+b-b_0)\frac{1}{N}\sum_{l\in L}e^{2\pi il\cdot(b-b_0)}.$$
If $(B,L)$ is a Hadamard pair then we get 
$A=f(x)$, which proves the other Cuntz relation 
$$\sum_{l\in L}S_lS_l^*=I_{L^2(\mu_B)}.$$
\end{proof}

\begin{definition}\label{def1.5}
Suppose $(B,L)$ form a Hadamard pair. We denote by $\Rs$, the transpose of the matrix $R$. We consider the {\it dual} IFS 
\begin{equation}
\tau_l^{(L)}(x)=(\Rs)^{-1}(x+l),\quad(x\in \br^d, l\in L)
\label{eq1.12}
\end{equation}
We denote the attractor of this IFS by $X_L$ and its invariant measure by $\mu_L$. 

To simplify the notation, we will use $\tau_l$ instead of $\tau_l^{(L)}$, and make the convention that when the subscript is an $l$ then we refer to $\tau_l^{(L)}$, and when the subscript is a $b$ we refer to $\tau_b$. 

\end{definition}

\begin{definition}\label{deftro}
For the set $B$ define the function 
\begin{equation}
\chi_B(x)=\frac{1}{N}\sum_{b\in B}e^{2\pi i b\cdot x},\quad(x\in \br^d)
\label{eq1.13}
\end{equation}

The {\it transfer operator} $R_{B,L}$ is defined on functions $f$ on $\br^d$ by 
\begin{equation}
(R_{B,L}f)(x)=\sum_{l\in L}|\chi_B(\tau_l (x))|^2f(\tau_l (x))=
\sum_{l\in L}|\chi_B(\Rsi(x+l))|^2f(\Rsi(x+l)).
\label{eqtro}
\end{equation}
\end{definition}

The operators $S_l^*$ behave well on exponential functions: the next lemma will be helpful for our computations.
\begin{lemma}\label{lem1.8}
The following assertions hold:
\begin{enumerate}
\item
For all $l\in L$ and $t\in \br^d$
\begin{equation}
S_l^*e_t=e_{-\tau_l^{(L)}(-t)}\chi_B(-\tau_l^{(L)}(-t)).
\label{eq1.8.1}
\end{equation}
\item
Let $\H$ be a reducing subspace of $L^2(\mu_B)$ for the representation $(S_l)_{l\in L}$. If $e_t\in\H$ and $\chi_B(-\tau_l^{(L)}(-t))\neq 0$, then $e_{-\tau_l^{(L)}(-t)}\in \H$.
\end{enumerate}

\end{lemma}

\begin{proof}
We have
$$(S_l^*e_t)(x)=\frac1N\sum_{b\in B}e^{-2\pi il\cdot R^{-1}(x+b)}e^{2\pi it\cdot R^{-1}(x+b)}=\frac1N\sum_{b\in B}e^{-2\pi i\Rsi(-t+l)\cdot(x+b)}$$
$$=e^{-2\pi i\Rsi(-t+l)\cdot x}\frac1N\sum_{b\in B}e^{-2\pi i\Rsi(-t+l)\cdot b}$$$$=e^{-2\pi i\tau_l^{(L)}(-t)\cdot x}\chi_B(-\tau_l^{(L)}(t))=e_{-\tau_l^{(L)}(x)}\chi_B(-\tau_l^{(L)}(-t)).$$

(ii) follows from (i).
\end{proof}

\begin{definition}\label{def1.9}
We will use the following notation for a finite word $w=w_1\dots w_n\in L^n$,
$$S_w=S_{w_1}\dots S_{w_n}.$$

\end{definition}

\begin{proposition}\label{prinvk}
For a subspace $H$ of $L^2(\mu_B)$ denote by $P_H$ the orthogonal projection onto $H$. 
Let $\K$ be a subspace which is invariant for all the maps $S_l^*$, $l\in L$, i.e., $S_l^*\K\subset \K$ for all $l\in L$. Let $\K_0:=\K$, 
$$\K_n:=\cj{\textup{span}}\left\{S_w\K : w\in L^n\right\},$$
$$\K_\infty:=\cj{\textup{span}}\left\{S_w\K : w\in L^n, n\in\bn\right\}.$$
Then 
\begin{enumerate}
	\item For all $n\in \bn$, $\K_n$ is invariant for $S_l^*$, $l\in L$ and $S_l^*\K_{n+1}=\K_n$ for all $l\in L$, $n\geq 0$.
	\item $P_{\K_n}\leq P_{\K_{n+1}}$ for all $n\in\bn$. 
	\item For all $n\in\bn$
	$$P_{\K_{n+1}}=\alpha(P_{\K_n})=\sum_{l\in L}S_lP_{\K_n}S_l^*.$$
	\item The projections $P_{\K_n}$ converge to $P_{\K_\infty}$ in the strong operator topology.
\end{enumerate}

\end{proposition}
\begin{proof}
(i) Take $w=w_1\dots w_n\in L^n$ and $k\in \K$. Then for any $l_0\in L$
$$S_{l_0}^*S_{w_1}S_{w_2\dots w_n}k=\delta_{l_0w_1}S_{w_2\dots w_n}k=\delta_{l_0 w_1}\sum_{l\in L}S_{w_2\dots w_n}S_lS_l^*k.$$
But $S_l^*k\in\K$ so $S_{l_0}^*S_wk\in\K_n$. 

 This computation implies also that $S_{l_0}^*\K_n\subset \K_{n-1}$. The other inclusion follows from the Cuntz relations, and this implies (i). 
 
 (ii) is immediate from (i).

 Since $\K_n$ is also invariant under the maps $S_l^*$, it is enough to prove (i) for $n=0$. Take $l_0\in L$ and $k\in\K$. Then 
 $$\alpha(P_{\K})S_{l_0}k =\sum_{l\in L}S_lP_{\K}S_l^*S_{l_0}k=S_{l_0}P_{\K}k=S_{l_0}k.$$
 Therefore $\K_1$ is contained in the range of the projection $\alpha(P_{\K})$. 
 
 Also, for any $\alpha(P_{\K})v=\sum_{l\in L}S_lP_{\K}S_l^*v\in\K_1$ since $P_{\K}S_l^*v\in\K$.

 (iv) is clear since $\cup_n\K_n$ spans $\K_\infty$.
\end{proof}

\begin{definition}\label{defht}
For an operator $T$ on $L^2(\mu_B)$ define the following function 
$$h_T(t)=\ip{Te_{-t}}{e_{-t}},\quad(t\in\br^d).$$
\end{definition}

\begin{proposition}\label{prrht}
Let $T$ be an operator on $L^2(\mu_B)$. Then 
	$$h_{\alpha(T)}=R_{B,L}h_T.$$
The function $h_T$ is entire analytic.

\end{proposition}

\begin{proof}
We have, with Lemma \ref{lem1.8}:
$$h_{\alpha(T)}(t)=\ip{\sum_{l\in L}S_lTS_l^*e_{-t}}{e_{-t}}=\sum_{l\in L}\ip{TS_l^*e_{-t}}{S_l^*e_{-t}}=\sum_{l\in L}|\chi_B(\tau_l (t))|^2\ip{Te_{-\tau_l (t)}}{e_{-\tau_l (t)}}$$
$$=\sum_{l\in L}|\chi_B(\tau_l (t))|^2h_T(\tau_l (t))=R_{B,L}h_T(t).$$

Since the operator $T$ is bounded, it is easy to check that the function $h_T$ is entire analytic. 
\end{proof}

\begin{corollary}\label{corcto}
Using the notations above:
\begin{enumerate}
\item
If $A$ is a bounded operator in the commutant of the representation $(S_l)_{l\in L}$ then the function 
$h_A$
is an entire analytic harmonic function for the transfer operator $R_{B,L}$, i.e.,
\begin{equation}
R_{B,L}h_A=h_A
\label{eqrha}
\end{equation}

\item

Suppose $\K$ is a subspace which is invariant under all the maps $S_l^*$, $l\in L$. Then the function $h_{P_{\K}}$ is an entire analytic subharmonic function for the transfer operator $R_{B,L}$ i.e., 
$$R_{B,L}h_{P_{\K}}\geq h_{P_{\K}}.$$
\end{enumerate}
\end{corollary}

\begin{proof}
If $A$ commutes with the representation then $\alpha(A)=\sum_l S_lAS_l^*=A\sum_lS_lS_l^*=A$. If $\K$ is invariant under the maps $S_l^*$ then 
$\alpha(P_{\K})\geq P_{\K}$ by Proposition \ref{prinvk}. Then the corollary follows directly from Proposition \ref{prrht}.

\end{proof}

\begin{corollary}\label{cor1.14}
Let $\K$ be a subspace of $L^2(\mu_B)$ which is invariant under all the maps $S_l^*$, $l\in L$. With the notation in Proposition \ref{prinvk}, the following limit exists uniformly on compact sets:
$$\lim_{n\rightarrow\infty}R_{B,L}^nh_{P_{\K}}=h_{P_{\K_\infty}}.$$
\end{corollary}

\begin{proof}
With Propositions \ref{prrht} and \ref{prinvk} we have
$$R_{B,L}^nh_{P_{\K}}=h_{\alpha^n(P_{\K})}=h_{P_{\K_n}}.$$
Since $P_{\K_n}$ converges in the strong operator topology to $P_{\K_\infty}$, we obtain that $h_{P_{\K_n}}$ converges to $h_{P_{\K_\infty}}$ pointwise. Using Corollary \ref{corcto}, we see that this is an increasing sequence of functions; then, with Dini's theorem we obtain the conclusion.
\end{proof}

\section{Representations associated to minimal words}

 We now turn to a dynamical systems feature which determines both the algorithmic and the analytic part of the problem of $L^2(\mu)$. This feature in turns divides up into two parts, periodic and non-periodic. The precise meaning of these terms is fleshed out in Definitions \ref{def2.1} and 
\ref{def2.2} (in the present section), and in Definition \ref{def4.2} inside the paper. In both cases, we deal with a random walk-dynamical system. The first case is especially easy to understand in terms of a natural encoding with finite and infinite code-words. By contrast, the second case involves invariant sets (for the walk, Definition \ref{def4.2}) which can have quite subtle fractal properties. Some of the possibilities (case 2) are illustrated by examples in section 4.

      We begin with a discussion of the minimal words in the encoding for case 1.

\begin{definition}\label{def2.1}
Consider a finite word over a finite alphabet $L$ with $\#L=N$. We say that $w$ is {\it minimal} if there is no word $u$ such that $w=\underbrace{uu\dots u}_{p\mbox{ times}}$ for some $p\geq 2$.

We denote by $\underline w=ww\dots$, the infinite word obtained by the repetition of $w$ infinitely many times. For two words $u$, $w$ with $u$ finite, we denote by $uw$ the concatenation of the two words. 

We define the shift $\sigma$ on infinite words $\sigma(\omega_1\omega_2\dots)=\omega_2\omega_3\dots$. 

\end{definition}

\begin{definition}\label{def2.2}
Let $w$ be a minimal word. Let $\Gamma(w)$ be the set of infinite words over the alphabet $L$ that end in an infinite repetition of the word $w$, i.e., 
\begin{equation}
\Gamma(w):=\{u\underline w : u\mbox{ is a finite word over }L\}
\label{eq2.1}
\end{equation}
We will use the Dirac notation for vectors in $l^2(\Gamma(w))$: for $\omega\in l^2(\Gamma(w))$
$$
\ket{\omega}=\delta_\omega,\quad \delta_\omega(\xi)=\left\{
\begin{array}{cc}
1,&\xi=\omega\\
0,&\xi\neq\omega
\end{array}\right.
$$

We define the represetation $\rho_w$ of the Cuntz algebra $\mathcal O_N$ by 
\begin{equation}
\rho_w(S_l)\ket{\omega}=\ket{l\omega},\quad (\omega\in l^2(\Gamma(w)))
\label{eq2.2}
\end{equation}
Moreover, if $u,v$ are vectors in a Hilbert space, we will use the notation $|u\rangle\langle v|$ for the corresponding rank-one operator.
\end{definition}

\begin{theorem}\label{th2.3} The following assertions hold:
\begin{enumerate}
\item
The operators in \eqref{eq2.2} defines an irreducible representation of the Cuntz algebra $\mathcal O_N$.
\item 
Define the subspace $\mathcal K_w$ spanned by the vectors $\{\ket{\underline w},\ket{\sigma(\underline w)},\dots, \ket{\sigma^{p-1}(\underline w)}\}$, where $p$ is  the length of $w$. Then $\mathcal K_w$ is invariant for the operators $\rho_w(S_l)^*$, $l\in L$ and it is cyclic for the representation $\rho_w$. 
\item 
If $w$ and $w'$ are two minimal words that are not a cyclic permutation of each other, then the representations $\rho_w$ and $\rho_{w'}$ are disjoint. 
\end{enumerate}
\end{theorem}

\begin{proof}
It is easy to check that 
$$\rho_w(S_l)^*\ket{\omega}=\left\{\begin{array}{cc}
0,& \mbox{ if }\omega_1\neq l\\
\ket{\sigma(\omega)},&\mbox{ if } \omega_1=l
\end{array}
\right.$$
for $\omega=\omega_1\omega_2\dots$.
From this it follows after a simple computation that $\rho_w$ is a representation of the Cuntz algebra.

Since $\sigma^p(\underline w)=\underline w$, it follows that $\rho_w(S_l)^*\mathcal K_w\subset \mathcal K_w$. Since every word in $\Gamma(w)$ ends in $\underline w$, it can be easily seen that $\mathcal K_w$ is cyclic for the representation.

It remains to check that the representation is irreducible and assertion (iii). For this we will use Theorem from \cite{BJKW01}, applied to our situation (see also \cite{BrJo97a, BrJo97b}):

\begin{theorem}\label{th2.4}
There is a bijective correspondence between 
\begin{enumerate}
\item
Operators $A$ that intertwine the representations $\rho_w$ and $\rho_{w'}$, i.e., 
$$\rho_{w'}(S_l)A=A\rho_w(S_l),\quad(l\in L).$$
\item
Fixed points of the map 
$$\Phi(C)=\sum_{l\in L}V_l'CV_l^*,\quad(C\in\mathcal B(\mathcal K_w,\mathcal K_{w'})),$$
where $V_l=P\rho_w(S_l)P$, $V_l'=P'\rho_{w'}(S_l)P'$, $l\in L$, with $P=$projection onto $\mathcal K_w$, $P'=$projection onto $\mathcal K_{w'}$. 
The correspondence from (i) to (ii) is given by 
$$C=P'AP.$$
\end{enumerate}

\end{theorem}

Returning to the proof of Theorem \ref{th2.3}, we will compute the fixed points of the map $\Phi$.

Let $u_j=\ket{\sigma^j(\underline w)}$, $j=0,\dots,p-1$, where $p=$length of $w$, and let $u_i'=\ket{\sigma^i(\underline w')}$, $i=0,\dots,p'-1$, where $p'=$length of $w'$. The space $\mathcal B(\mathcal K_w,\mathcal K_{w'})$ is spanned by the rank one operators $\ket{u_i'}\bra{u_j}$. We have
$$\Phi(\ket{u_i'}\bra{u_j})=\sum_{l\in L}\ket{V_l'u_i'}\bra{V_lu_j}$$
But 
$$V_lu_j=P\rho_w(S_l)u_j=P\ket{l\sigma^j(\underline w)}=P\ket{lw_{j+1}w_{j+2}\dots w_p\underline w}=\delta_{l,w_j}\ket{\sigma^{j-1}(\underline w)}$$ 
(we use here a notation $\mod p$, or $\mod p'$ when required, so $\sigma^{-1}(\underline w)$ will mean $\sigma^{p-1}(\underline w)$)

Then 
$$\Phi(\ket{u_i'}\bra{u_j})=\sum_{l\in L}\delta_{l,w_i'}\delta_{l,w_j}\ket{\sigma^{i-1}(\underline w')}\bra{\sigma^{j-1}(\underline w)}=\delta_{w_i',w_j}\ket{\sigma^{i-1}(\underline w')}\bra{\sigma^{j-1}(\underline w)}.$$

Now suppose $C=\sum_{i=0}^{p'-1}\sum_{j=0}^{p-1}c_{i,j}\ket{u_i'}\bra{u_j}$ is a fixed point for $\Phi$. Then 
$$\sum_{i=0}^{p'-1}\sum_{j=0}^{p-1}c_{i,j}\ket{u_i'}\bra{u_j}=\sum_{i,j}c_{i,j}\delta_{w_i',w_j}\ket{u_{i-1}'}\bra{u_{j-1}}.$$
Then 
$$c_{i,j}=c_{i+1,j+1}\delta_{w_{i+1}',w_{j+1}},\mbox{ for all }i,j.$$
Thus, $c_{i,j}\neq 0$ only if $w_{i+k}'=w_{j+k}$ for all $k$, which means that $\sigma^{i-j}(\underline w')=\underline w$. But this implies, since $w,w'$ are minimal, that $w$ is a cyclic permutation of $w'$. 
Thus, if $w$, $w'$ are not cyclic permutations of each other, then the only fixed point of $\Phi$ is $C$, and with Theorem \ref{th2.4}, this implies that there are no nonzero intertwining operators.

To check that the representation $\rho_w$ is irreducible, we use the same computation, now with $w=w'$. We saw that $c_{i,j}\neq 0$ only if $\sigma^{i-j}(\underline w')=\underline w$. But since $|i-j|<p$ and $w$ is minimal this implies that $i=j$. In the case $i=j$, the same computation implies that $c_{i,i}=c_{i+k,i+k}$, therefore $C$ is a scalar multiple of the identity. Using again Theorem \ref{th2.4}, we obtain that the only operators in the commutant of $\rho_w$ are scalar multiples of the identity.

\end{proof}

\section{Representations associated to extreme cycles}

 We mentioned in section 2, that the algorithmic and the analytic part of the problem of $L^2(\mu)$ involves a dynamical system. Its nature divides up into two parts, periodic and non-periodic. Case 1 has an algorithmic part (code-words), and an analytic part taking the form of extreme cycles, and their associated representations. This is worked out below.

\begin{definition}\label{defcyc}

A finite set of distinct points $C=\{x_0,x_1,\dots,x_{p-1}\}$ is called an $L$-cycle if there exist $l_0,\dots, l_{p-1}\in L$ such that 
$$\tau_{l_0} (x_0)=x_1, \tau_{l_1} (x_1)=x_2,\dots, \tau_{l_{p-2}} (x_{p-2})=x_{p-1}, \tau_{l_{p-1}} (x_{p-1})=x_{0}.$$
We call $w(C):=l_0l_1\dots l_{p-1}$ the {\it word of the cycle }$C$.
The points in $C$ are called $L$-cycle points. 

An $L$-cycle is called {\it $B$-extreme} if 
\begin{equation}
\left|\chi_B(x)\right|=1,\mbox{ for all } x\in C.
\label{eq1.14}
\end{equation}
\end{definition}

\begin{theorem}\label{th1.6}\cite[Theorem 4.1]{DJ07a}
Under the conditions above, assume $(B,L)$ is a Hadamard pair. Suppose there exist $d$ linearly independent vectors in the set 
\begin{equation}
\Gamma(B):=\left\{\sum_{k=0}^nR^kb_k : b_k\in B, n\in \bn\right\}.
\label{eq1.15}
\end{equation}
Define 
\begin{equation}
\Gamma(B)^\circ:=\left\{ x\in\br^d : \beta\cdot x\in\bz\mbox{ for all }\beta\in\Gamma(B)\right\}.
\label{eq1.16}
\end{equation}
Then $\Gamma(B)^\circ$ is a lattice that contains $\bz^d$, is invariant under $\Rs$, and if $l,l'\in L$ with $l-l'\in \Rs\Gamma(B)^\circ$ then $l=l'$. Moreover 
\begin{equation}
\Gamma(B)^\circ\cap X_L=\bigcup\left\{ C : C\mbox{ is a $B$-extreme $L$-cycle}\right\}.
\label{eq1.17}
\end{equation}

\end{theorem}

\begin{definition}\label{def1.7}
Let $(B,L)$ be a Hadamard pair. We say that the Hadamard pair is {\it regular} if the IFS $(\tau_b)_{b\in B}$ has no overlap and there exist $d$ linearly independent vectors in the set $\Gamma(B)$ defined by \eqref{eq1.15}.

Define the maps $\sigma_l$ on $\br^d$ by
\begin{equation}
\sigma_l(x)=\Rs x+l,\quad (x\in\br^d,l\in L)
\label{eq1.18}
\end{equation}

A subset $\Lambda$ of $\br^d$ is called {\it invariant} if $\sigma_l(\Lambda)\subset\Lambda$ for all $l\in L$. 
\end{definition}

\begin{theorem}\label{th1.8}
Let $(B,L)$ be a regular Hadamard pair and let $C$ be a $B$-extreme $L$-cycle. Let $\Lambda(C)$ be the smallest invariant set that contains $-C$. Define the subspace of $L^2(\mu_B)$:
\begin{equation}
\H(C):=\cj{\operatorname*{span}}\left\{e_\lambda : \lambda\in\Lambda(C)\right\}.
\label{eq1.19}
\end{equation}
Then
\begin{enumerate}
\item 
The space $\H(C)$ is a reducing subspace for the representation $(S_l)_{l\in L}$ of the Cuntz algebra $\mathcal O_N$.
\item 
The exponential function $\{e_\lambda : \lambda\in\Lambda(C)\}$ form an orthonormal basis for $\H(C)$ and 
\begin{equation}
S_l(e_\lambda)=e_{\sigma_l(\lambda)},\quad(\lambda\in\Lambda(C),l\in L)
\label{eq1.19.1}
\end{equation}
\item
If $C_1$ and $C_2$ are two distinct $B$-extreme $L$-cycles, then the subspaces $\H(C_1)$ and $\H(C_2)$ are orthogonal.
\item 
Let $w(C)$ be the word of the cycle $C$. The restriction of the representation $(S_l)_{l\in L}$ to $\H(C)$ is equivalent to the representation $\rho_{w(C)}$ on $l^2(\Gamma(w(C))$ from Theorem \ref{th2.3}. The isomorphism can be defined as follows: define the encoding map $\epsilon_C:\Lambda(C)\rightarrow\Gamma(w(C))$, $\epsilon_C(\lambda)=\omega_1\omega_2\dots$, where $\omega_1,\omega_2,\dots $ are uniquely defined by the condition 
$$S_{\omega_n}^*\dots S_{\omega_1}^*e_\lambda\neq 0,\quad(n\in\bn).$$
Then the map 
$W_C:\H(C)\rightarrow l^2(\Gamma(w(C))$
$$W_C(e_\lambda)=\ket{\epsilon_C(\lambda)}$$
defines an isometric isomorphism that intertwines the two representations of $\mathcal O_N$.

\item 
The restrictions of the representation $(S_l)_{l\in L}$ to the subspaces $\H(C)$, where $C$ is any $B$-extreme $L$-cycle, are mutually disjoint irreducible representations of the Cuntz algebra $\mathcal O_N$. 
\end{enumerate}

\end{theorem}

We begin with some lemmas:
\begin{lemma}\label{lem1.9}
In the hypotheses of Theorem \ref{th1.8} we have:
\begin{enumerate}
\item
$\Lambda(C)\subset \Gamma(B)^\circ$.
\item For every $t\in \Gamma(B)^\circ$ and $l\in L$, $S_le_t=e_{\sigma_l(t)}$.
\item If $t=\sigma_l(t')$ for some $l\in L$ and $t'\in \Gamma(B)^\circ$, then $S_{l'}^*e_t=\delta_{ll'}e_{t'}$.
\item For all $t\in\Lambda(C)$, there exist a unique $t'\in\Lambda(C)$ and $l\in L$ such that $t=\sigma_l(t')$.
Moreover, if $t\in -C$ then $t'\in -C$.
\end{enumerate}
\end{lemma}

\begin{proof}
(i) follows from Theorem \ref{th1.6}. For (ii), let $t\in\Gamma(B)^\circ$ and $x\in\tau_b(X_B)$ for some $b\in B$. Then
$$(S_le_t)(x)=e^{2\pi il\cdot x}e^{2\pi it\cdot \R x}=e^{2\pi i l\cdot x}e^{2\pi it\cdot(\Rs x-b)}$$
But $t\in \Gamma(B)^0$ implies that $t\cdot b\in\bz$ so 
$$(S_le_t)(x)=e_{\sigma_l(t)}(x).$$

(iii) If $t=\sigma_l(t')$ with $t'\in\Gamma(B)^\circ$ then, using (ii),
$$S_{l'}^*e_t=S_{l'}^*S_le_{t'}=\delta_{ll'}e_{t'}.$$

(iv) It is clear that 
$$\Lambda(C)=\left\{\sigma_{l_n}\dots\sigma_{l_1}x : x\in -C, l_1,\dots,l_n\in L, n\in\bn\right\}.$$
Therefore, for existence, we have to check only that if $x\in -C$ then $x=\sigma_l(t')$ for some $l\in L$ and $t'\in -C$. Since $-x\in C$ and $C$ is an $L$-cycle, there exist $l\in L$ and $y\in C$ such that $\Rsi(x+l)=y$. This implies that $\sigma_l(-y)=-x$, and this proves the assertion. 

For uniqueness, suppose $t=\sigma_{l_1}(t_1')=\sigma_{l_2}(t_2')$ with $l_1,l_2\in L$, $t_1,t_2\in\Lambda(C)$. Then 
$l_1-l_2=\Rs(t_2'-t_1')$. But since $t_2'-t_1'\in \Gamma(B)^\circ$, Theorem \ref{th1.6} implies that $l_1=l_2$ and so $t_1'=t_2'$.
\end{proof}

\begin{lemma}\label{lem1.10}
Let $(B,L)$ be a regular Hadamard pair. Then 
\begin{enumerate}
\item The Fourier transform of the measure $\mu_B$, defined by 
\begin{equation}
\widehat\mu_B(x)=\int e^{2\pi i t\cdot x}\,d\mu_B(t),\quad(t\in\br^d)
\label{eq1.20}
\end{equation}
satisfies the following equality
\begin{equation}
\widehat\mu_B(x)=\prod_{n=1}^\infty\chi_B((\Rs)^{-n}x),\quad(x\in \br^d).
\label{eq1.21}
\end{equation}
\item Two exponential functions $e_t$ and $e_{t'}$ are orthogonal in $L^2(\mu_B)$ iff $\widehat\mu_B(t-t')=0$.
\item The function $\chi_B$ satisfies the following {\it QMF condition}
\begin{equation}
\sum_{l\in L}|\chi_B(\Rsi(x+l))|^2=1,\quad(x\in\br^d).
\label{eq1.22}
\end{equation}
\item For all $t\in\Gamma(B)^\circ$
\begin{equation}
\chi_B(x+t)=\chi_B(x),\quad(x\in\br^d)
\label{eq1.23}
\end{equation}
\item If $C$ is a $B$-extreme $L$-cycle, and $x\in C$ with $-x=\sigma_{l_0}(-y)$ for some $l_0\in L$ and $y\in C$, then 
\begin{equation}
\chi_B(\Rsi(x+l))=0\mbox{ for all }l\in L, l\neq l_0.
\label{eq1.24}
\end{equation}
\end{enumerate}

\end{lemma}

\begin{proof}
For the proof of (i)--(iii) we refer to \cite{DJ06b}. For (iv), if $t\in \Gamma(B)^\circ$ then $b\cdot t\in\bz$ so 
$e^{2\pi ib\cdot t}=1$ for all $b\in B$ and this implies (iv). 

For (v), the equality $-x=\sigma_{l_0}(-y)$ implies that $\Rsi(x+l_0)=y$. Since the cycle is $B$-extreme $|\chi_B(y)|=1$. Then, using (iii), one obtains that $|\chi_B(\Rsi(x+l)|=0$ for all $l\neq l_0$.

\end{proof}

\begin{proof}[Proof of Theorem \ref{th1.8}]
The fact that $\H(C)$ is invariant for all $S_l$ and $S_l^*$ follows from Lemma \ref{lem1.9}(ii) and (iii).

Next let $C'$ be a $B$-extreme $L$-cycle, and let $t\in\Lambda(C)$, $t'\in\Lambda(C')$. Then there exist $x\in C$, $x'\in C'$, $l_1,\dots, l_n\in L, l_1',\dots,l_m'\in L$ such that $t=\sigma_{l_n}\dots\sigma_{l_1}(-x)$ and $t'=\sigma_{l_n'}\dots\sigma_{l_1'}(-x')$. Using Lemma \ref{lem1.9}(iv), and composing with a few more $\sigma_l$'s if necessary, we can assume $n=m$. We have then with Lemma \ref{lem1.9}(ii), and then using the Cuntz relations, 
$$\ip{e_t}{e_{t'}}=\ip{S_{l_n}\dots S_{l_1}e_{-x} }{S_{l_n'}\dots S_{l_1'}e_{-x'}}=\delta_{l_nl_n'}\dots\delta_{l_1l_1'}\ip{e_{-x}}{e_{-x'}}=\delta_{l_nl_n'}\dots\delta_{l_1l_1'}\widehat\mu_B(x'-x).$$
Thus, to prove (ii) and (iii), it is enough to show that $\widehat\mu_B(x'-x)=0$ for $x\neq x'$. 

By Lemma \ref{lem1.9}(iv), $-x=\sigma_{l_0}(-y)$ for some $l_0\in L$ and $y\in C$, and $-x'=\sigma_{l_0'}(-y')$ for some $l_0'\in L$ and $y'\in C'$. We can assume $l_0\neq l_0'$. If not, then one can compose with a few more $\sigma_l$'s in the argument above until the cycle $C$ uses another digit $l_0$ than the cycle $C'$. 

We have 
$$\chi_B(\Rsi(x'-x))=\chi_B\left(\Rsi(x'+l_0)-\Rsi(x+l_0)\right)=\chi_B\left(\Rsi(x'+l_0)-y\right)$$
$$=\chi_B(\Rsi(x'+l_0))\quad\mbox{ (because $y\in\Gamma(B)^\circ$ and due to Lemma \ref{lem1.10}(iv))}$$
$$=0,\quad \mbox{(from Lemma \ref{lem1.10}(v), because $l_0\neq l_0'$)}.$$

To check (iv), we note that 
$$W_CS_le_\lambda=W_Ce_{\sigma_l(\lambda)}=\ket{\epsilon_C(\sigma_l(\lambda))}=\ket{l\epsilon_C(\lambda)}=\rho_{w(C)}(S_l)\ket{\epsilon_C(\lambda)}=\rho_{w(C)}(S_l)W_Ce_\lambda.$$
Thus $W_C$ is intertwining.

Also, we have $\tau_{w(C)_1}x_0=x_1$ so $-x_0=\Rs (-x_1)+w(C)_1$, and therefore $S_{w(C)_1}e_{-x_1}=e_{-x_0}$ so $S_{w(C)_1}^*e_{-x_0}=e_{-x_1}$. By induction, we get that
$\epsilon_C(-x_0)=\underline{w(C)}$.

From the definition of $\Lambda(C)$ we see that for each $\lambda\in \Lambda(C)$, there exists a point $x_k$ in $-C$ and some digits $l_1,\dots, l_n$ such that 
$\lambda=\sigma_{l_1}\dots\sigma_{l_n}(-x_k)$. Composing with a few more $\sigma_l$'s we can assume $x_k=x_0$, so $\lambda =\sigma_{l_1}\dots\sigma_{l_n}(-x_0)$.

Then $$S_{l_1}e_{\sigma_{l_{2}}\dots \sigma_{l_n}(-x_0)}=e_\lambda,\mbox{ so } S_{l_1}^*e_\lambda=e_{\sigma_{l_{2}}\dots \sigma_{l_n}(-x_0)}.$$
By induction, we get that $\epsilon_C(\lambda)=l_1\dots l_{n}\underline{w(C)}$.

This proves in particular that $\epsilon_C$ is surjective, so $W_C$ is onto. 

Note that we cannot have $\lambda=\sigma_{l_1}\lambda_1=\sigma_{l_1'}\lambda_1'$ for $l_1\neq l_1'$ and $\lambda_1,\lambda_1'\in\Lambda(C)$, because in this case 
$l_1+\Rs\lambda_1=l_1'\Rs\lambda_1'$ so $l_1-l_1'\in \Rs\Gamma(B)^\circ$, and with Theorem \ref{th1.6}, it follows that $l_1=l_1'$. Therefore, the encoding $\epsilon_C$ is well defined.

To see that $\epsilon_C$ is injective, suppose $\epsilon_C(\lambda)=\epsilon_C(\lambda')=l_1\dots l_n\underline{w(C)}$ then $\lambda=\sigma_{l_1}\dots\sigma_{l_n}(-x_0)=\lambda'$. 
Hence $W_C$ maps an orthonormal basis to an orthobormal basis and therefore it is an isometric isomorphism.

The last statement in the theorem follows now from Theorem \ref{th2.3}.

\end{proof}

\begin{corollary}\label{cor3.7}
If the dimension $d=1$ then the decomposition of the representation $(S_l)_{l\in L}$ into irreducible subrepresentation is given by 
$$L^2(\mu_B)=\bigoplus \left\{\H(C) : C\mbox{ is a $B$-extreme $L$-cycle}\right\}.$$
The commutant is finite-dimensional and abelian.

\end{corollary}

\begin{proof}
From \cite{DJ06b} we know that, when the dimension $d=1$, then the union of the sets $\Lambda(C)$ is a spectrum for $\mu_B$, therefore the corresponding exponentials form a complete orthonormal set in $L^2(\mu_B)$. Everything follows then from Theorem \ref{th1.8}.
\end{proof}

\section{Beyond cycles}

While in simple cases, the harmonic analysis of $L^2(\mu)$ may be accounted for by cycles (sections 2 and 3), there is a wider class involving more non-cyclic invariant sets (Definition \ref{def4.2}). These are studied below. We show (Theorem \ref{pr4.4}) that even in the non-cyclic case we still get an associated splitting of the Hilbert space $L^2(\mu)$ into orthogonal closed subspaces. We further show by examples that the invariant sets can have quite subtle fractal properties. Some of the possibilities (case 2) are illustrated in Example \ref{ex4.7}, and the results following it. This approach to invariant sets was first introduced in \cite{CoRa90,CCR96,CHR97}.
\begin{definition}\label{def4.1}
For a finite word $w=l_1\dots l_n$ we use the notation
$$\tau_w=\tau_{l_n}\dots\tau_{l_1}.$$
$$\chi_B^{(n)}(x)=\chi_B(x)\chi_B(\Rs x)\dots\chi_B((\Rs)^{n-1}x),\quad(x\in\br^d)$$

For the operators $(S_l)_{l\in L}$ defined in \eqref{eq1.9} and any word $w=l_1\dots l_n$ we use the notation
$$S_w=S_{l_1}\dots S_{l_n}.$$
\end{definition}

\begin{definition}\label{def4.2}
A subset $M$ of $\br^d$ is called $L$-invariant if for all $x\in M$ and all $l\in L$, if $\chi_B(\tau_l(x))\neq0$ then $\tau_l x\in M$. We say that the transition from $x$ to $\tau_lx$ is possible if $\chi_B(\tau_lx)\neq0$. For $x\in \br^d$ we denote by $O(x)$, the orbit of $x$, that is the set of all points $y\in \br^d$ such that there exist $y_1=x_0,\dots, y_n=y$ such that the transition from $y_i$ to $y_{i+1}$ is possible.
Two $L$-invariant sets $M_1,M_2$ are called separated if $\dist(M_1,M_2)>0$. (Here $\dist$ denotes a distance under which all maps $\tau_l$ are strict contractions, and which generates the Euclidian topology on $\br^d$.)
A family $\{M_i : i=1,\dots,n\}$ of closed $L$-invariant sets is called complete if for any closed $L$-invariant set $M$, 
$$M\cap\bigcup_{i=1}^n M_i\neq \emptyset.$$

A closed $L$-invariant set $M$ is called minimal if for any closed $L$-invariant set $M'\subset M$, $M'=M$. 
\end{definition}

Here are a few simple properties of $L$-invariant sets.

\begin{proposition}\label{pr4.2.1}
The following assertions hold:
\begin{enumerate}
\item If $M$ is an $L$-invariant set, then the closure $\cj M$ is also $L$-invariant. 
\item Any $B$-extreme $L$-cycle is an $L$-invariant set. 
\item For any minimal $L$-invariant set $M$ and any $x\in M$, $M=\cj{O(x)}$. Also $M$ is contained in the attractor $X_L$ of the IFS $(\tau_l)_{l\in L}$.
\item The family of minimal closed $L$-invariant sets is a finite, complete family of separated $L$-invariant sets. 

\end{enumerate}

\end{proposition}

\begin{proof}
(i) Take $x\in \cj M$ and $l\in L$ such that $\chi_B(\tau_lx)\neq 0$. There exists a sequence in $x_n\in M$ that converges to $x$. For $n$ large enough $\chi_B(\tau_lx_n)\neq 0$. Since $M$ is invariant, it follows that $\tau_lx_n\in M$ so $\tau_lx\in\cj M$. Therefore $\cj M$ is invariant.

(ii) Since $(B,L)$ form a Hadamard pair, it follows (see \cite{DJ06b}) that we have the equality
\begin{equation}
\sum_{l\in L}|\chi_B(\tau_lx)|^2=1,\quad(x\in\br^d).
\label{eq4.2.1}
\end{equation} 
Let $C=\{x_0,x_1,\dots, x_{p-1}\}$ be a $B$-extreme $L$-cycle, with $\tau_{l_0}x_0=x_1,\dots, \tau_{l_{p-1}}x_{p-1}=x_0$. 
Since $|\chi_B(x_1)|=1$, it follows that $\chi_B(\tau_lx_0)\neq 0$ for $l\neq 0$. So the only possible transitions is from $x_0$ to $x_1$. Similarly for the other points of the cycle. Therefore $C$ is $L$-invariant.

(iii) Clearly $O(x)$ is an $L$-invariant set. Since $M$ is minimal and $M$ contains $O(x)$ it follows that $M=\cj{O(x)}$. Also $X_L$ is clearly $L$-invariant. Therefore $M\cap X_L$ is $L$-invariant. So $M\cap X_L$ is either $M$ or $\emptyset$. We will prove that it cannot be $\emptyset$.

Take $x\in M$. Using \eqref{eq4.2.1}, we see that a transition is always possible so we can construct $l_1,l_2,\dots $ such that 
$y_n:=\tau_{l_n}\dots\tau_{l_1}x\in O(x)$ for all $n$. But since $X_L$ is the attractor of the IFS $(\tau_l)_{l\in L}$ it follows that $\dist(y_n, X_L)\rightarrow 0$. So $\cj O(x)\cap X_L\neq \empty$ and therefore $M\cap X_L\neq\emptyset$.

(iv) First we prove that any two distinct minimal closed $L$-invariant sets $M_1,M_2$ are separated. Since $M_1,M_2$ are contained in $X_L$, they are compact, so it is enough to show that they are disjoint. But $M_1\cap M_2$ is also $L$-invariant and closed so $M_1\cap M_2=\emptyset$ due to minimality.

Next we prove that actually there exists a $\delta>0$ that does not depend on $M_1,M_2$ such that $\dist(M_1,M_2)\geq \delta$. 

Since $\chi_B$ is uniformly continuous on $X_L$, there exists a $\delta>0$ such that if $\dist(x,y)<\delta$ then $|\chi_B(x)-\chi_B(y)|<1/2N$. 

Suppose there exist two points $x_0\in M_1$, $y_0\in M_2$ such that $\dist(x,y)<\delta$. From \eqref{eq4.2.1}, there exists $l_0\in L$ such that $|\chi_B(\tau_{l_0}x_0)|\geq 1/N$. Then, since $\dist(\tau_{l_0}x_0,\tau_{l_0}y_0)<\dist(x,y)<\delta$ it follows that $|\chi_B(\tau_{l_0}y_0)|>1/2N$. So the transition from $y_0$ to $\tau_{l_0}y_0$. By induction, we can find $l_1,\dots, l_n,\dots$ such that the transitions from $\tau_{l_n}\dots\tau_{l_0}x_0$ to $\tau_{l_{n+1}}\dots\tau_{l_0}x_0$ and from $\tau_{l_n}\dots\tau_{l_0}y_0$ to $\tau_{l_{n+1}}\dots\tau_{l_0}y_0$ are possible. But 
$$\dist(\tau_{l_n}\dots\tau_{l_0}x_0,\tau_{l_n}\dots\tau_{l_0}y_0)\rightarrow 0.$$
Since $M_1, M_2$ are invariant this implies that $\dist(M_1,M_2)=0$, a contradiction. 
 
Thus $\dist(M_1,M_2)\geq \delta$. Since all minimal closed $L$-invariant sets are contained in the compact set $X_L$, it follows that there are only finitely many of them.

Finally, we have to prove that the family of minimal sets is complete. For this take a closed invariant set $M$. Using Zorn's lemma, $M$ contains a minimal closed $L$-invariant set. This proves the completeness of the family. 
\end{proof}

\begin{lemma}\label{lem4.3}
Let $M_1,M_2$ be two $L$-invariant sets with $\dist(M_1,M_2)>0$. Then for all $t_1\in M_1$ and $t_2\in M_2$ there exists $n\in\bn$ such that for all words $w\in L^n$ of length $n$
$$\chi_B^{(n)}(\tau_wt_1)\chi_B^{(n)}(\tau_wt_2)=0.$$
\end{lemma}

\begin{proof}
Let $d_0:=\dist(M_1,M_2)$. We can assume that $\Rsi$ is contractive under the metric $\dist$. Pick $n$ such that $c^n \dist(t_1,t_2)< d_0$, where $0<c<1$ is a Lipschitz constant for all maps $\tau_l$, i.e.,
$$\dist(\tau_lx,\tau_ly)\leq c\dist(x,y),\quad(x,y\in\br^d).$$

Consider a word $w=l_1\dots l_n$ of length $n$. Suppose 
$$\chi_B^{(n)}(\tau_wt_1)\chi_B^{(n)}(\tau_wt_2)\neq0.$$
Since $R$ and $L$ have integer entries, it follows that for $p=0,\dots,n-1$,
$$\chi_B((\Rs)^p\tau_wt_1)=\chi_B(\tau_{l_{n-p}}\dots\tau_{l_1}t_1).$$
Then we have for $p=0,\dots,n-1$
$$\chi_B(\tau_{l_{n-p}}\dots\tau_{l_1}t_i)\neq0,\quad(i=1,2).$$
Since $M_1,M_2$ are invariant, this implies that $\tau_wt_1\in M_1$, $\tau_wt_2\in M_2$. But then 
$$d_0\leq \dist(\tau_wt_1,\tau_wt_2)\leq c^n\dist(t_1,t_2)<d_0.$$
This contradiction implies our lemma.
\end{proof}

\begin{definition}\label{def4.4}
Let $M$ be an $L$-invariant set. Define 
$$\mathcal K_M:=\clspan \{ e_{-t} : t\in M\}\subset L^2(\mu_B).$$
$$\H(M):=\clspan\{ S_w e_{-t} : t\in M, w\in L^n, n\in \bn\}.$$

\end{definition}

\begin{theorem}\label{pr4.4} 

Let $M$ be a closed $L$-invariant set. 
\begin{enumerate}
\item The space $\K_M$ is invariant under all the maps $S_l^*$, $l\in L$.
\item The space $\H(M)$ is reducing for the representation $(S_l)_{l\in L}$ of the Cuntz algebra $\mathcal O_N$. 
\item
Let $M_1$ and $M_2$ be two closed $L$-invariant sets with $\dist(M_1,M_2)>0$. Then the subspaces $\H(M_1)$ and $\H(M_2)$ are orthogonal. 
\item Let $(M_i)_{i=1}^n$ be a complete family of closed separated $L$-invariant sets. Then
$$\bigoplus_{i=1}^n\H(M_i)=L^2(\mu_B).$$  
\end{enumerate}
\end{theorem}

\begin{proof}
From Lemma \ref{lem1.8}, we have 
$$S_l^*e_{-t}=\chi_B(-\tau_l(t))e_{-\tau_l(t)}=\cj{\chi_B(\tau_l(t))}e_{-\tau_l(t)}$$ 
for $t\in M$ and $l\in L$. So, if $\chi_B(\tau_l(t))=0$ then $S_l^*e_{-t}=0$, otherwise, since $M$ is invariant $\tau_l(t)\in M$ so $S_l^*e_{-t}\in \K_M$. 
This proves (i).

(ii) follows from (i) and the Cuntz relations. 

We prove now (iii). Iterating Lemma \ref{lem1.8}, we have for a word $w\in L^n$ and $t\in\br^d$:
\begin{equation}
S_w^*e_{-t}=\chi_B^{(n)}(-\tau_wt)e_{-\tau_wt}.
\label{eq4.1}
\end{equation}

First we prove that the subspaces $\K_{M_1}$ and $\K_{M_2}$ are orthogonal. Let $t_1\in M_1$ and $t_2\in M_2$. Take $n$ as in Lemma \ref{lem4.3}. For the $L^2(\mu_B)$-inner product we have:
$$\ip{e_{-t_1}}{e_{-t_2}}=\ip{e_{-t_1}}{\sum_{w\in L^n}S_wS_w^*e_{-t_2}}=\sum_{w\in L^n}\ip{S_w^*e_{-t_1}}{S_w^*e_{-t_2}}$$
$$=\sum_{w\in L^n}\chi_B^{(n)}(-\tau_wt_1)\cj{\chi_B^{(n)}(-\tau_wt_2)}\ip{e_{-\tau_wt_1}}{e_{-\tau_wt_2}}=0.$$

So the spaces $\K_{M_1}$ and $\K_{M_2}$ are orthogonal. Next, we prove that the spaces $\H(M_1)$ and $\H(M_2)$ are orthogonal. Take $w_1$ and $w_2$ two words over $L$. We can assume the length of $w_1$ is bigger than the length of $w_2$. Take $t_1\in M_1$ and $t_2\in M_2$. 

We have from the Cuntz relations: $S_{w_2}^*S_{w_1}=0$ or $S_{w_2}^*S_{w_1}=S_w$ where $w$ is a subword of $w_1$. Only the second case requires some computations.

With \eqref{eq4.1} we have:

$$\ip{S_{w_1}e_{-t_1}}{S_{w_2}e_{-t_2}}=\ip{S_w^*e_{-t_1}}{e_{-t_2}}=\chi_B^{(n)}(-\tau_wt_1)\ip{e_{-\tau_wt_1}}{e_{-t_2}}.$$
So, either $\chi_B^{(n)}(-\tau_wt_1)\neq 0$, in which case , since $M_1$ is invariant, $\tau_wt_1\in M_1$ so $e_{-\tau_wt_1}$ is in $\K_{M_1}$ and therefore the result of the previous computation is 0, or, $\chi_B^{(n)}(-\tau_wt_1)=0$ in which case the result is again 0. 

Thus, (iii) follows.

We prove now (iv). Let $P$ be the projection onto the orthogonal sum of the subspaces $\H(M_i)$. Since the subspaces are reducing, we have that $P$ is in the commutant of the representation. With Corollary \ref{corcto}, the function 
$$h_P(t)=\ip{Pe_{-t}}{e_{-t}}_{L^2(\mu_B)},\quad(t\in \br^d),$$
is a fixed point of the transfer operator $R_{B,L}$. Assume $h_P$ is not identically 1. We clearly have $0\leq h_P\leq 1$.

Take $r>c\max_{l\in L}\|l\|/(1-c)$, where $\|\cdot\|$ is a norm under which the maps $\tau_l$ are contractive, and $c$ is the contractive constant. Then for the closed ball $\cj B(0,r)$, one has 
$$\bigcup_{l\in L}\tau_l\cj B(0,r)\subset \cj B(0,r).$$

Let $m:=\min_{t\in \cj B(0,r)} h_P(t)$. Consider the set 
$$M:=\{ t\in \cj B(0,r) : h_P(t)=m\}.$$
Take $t\in M$. We have
$$m=h_P(t)=(R_{B,L}h_P)(t)=\sum_{l\in L}|\chi_B(\tau_lt)|^2h_P(\tau_lt)\geq m\sum_{l\in L}|\chi_B(\tau_lt)|^2=m,$$
it follows that, if $\chi_B(\tau_lt)\neq 0$ then $h_P(\tau_lt)=m$. This shows that $M$ is a closed $L$-invariant set. But, since the family $(M_i)_i$ is complete, it follows that there is a $t\in M\cap M_i$ for some $i\in \{1,\dots,n\}$. 

But, if $t\in M_i$ then $e_{-t}\in \H(M_i)$ so $h_P(t)=1$. On the other hand $t\in M$ so $h_P(t)=m$. So $m=1$. Therefore $h_P$ is constant $1$ on $\cj B(0,r)$. For $r\rightarrow\infty$, we get that $h_P=1$ on $\br^d$. This implies that $e_{-t}$ is in the range of $P$ for all $t\in\br^d$. But, by the Stone-Weierstrass theorem $e_{-t}$, $t\in\br^d$ span the entire space $L^2(\mu_B)$. So (iv) follows.  
\end{proof}

\begin{example}\label{ex4.7}
Consider the following data:
$$R=\begin{bmatrix}
2&0\\
1&2\end{bmatrix},\quad
B=\left\{ \begin{pmatrix}
0\\0
\end{pmatrix},
\begin{pmatrix}
1\\0
\end{pmatrix}
,
\begin{pmatrix}
0\\p
\end{pmatrix}
,
\begin{pmatrix}
1\\p
\end{pmatrix}
\right\}\quad
L=\left\{
\begin{pmatrix}
0\\0
\end{pmatrix}
,
\begin{pmatrix}
1\\0
\end{pmatrix}
,
\begin{pmatrix}
0\\1
\end{pmatrix}
,
\begin{pmatrix}
1\\1
\end{pmatrix}
\right\}$$
where $p\in\bn$ is an odd number.

\begin{figure}[ht]\label{fig1}
\centerline{ \vbox{\hbox{\epsfxsize 5cm\epsfbox{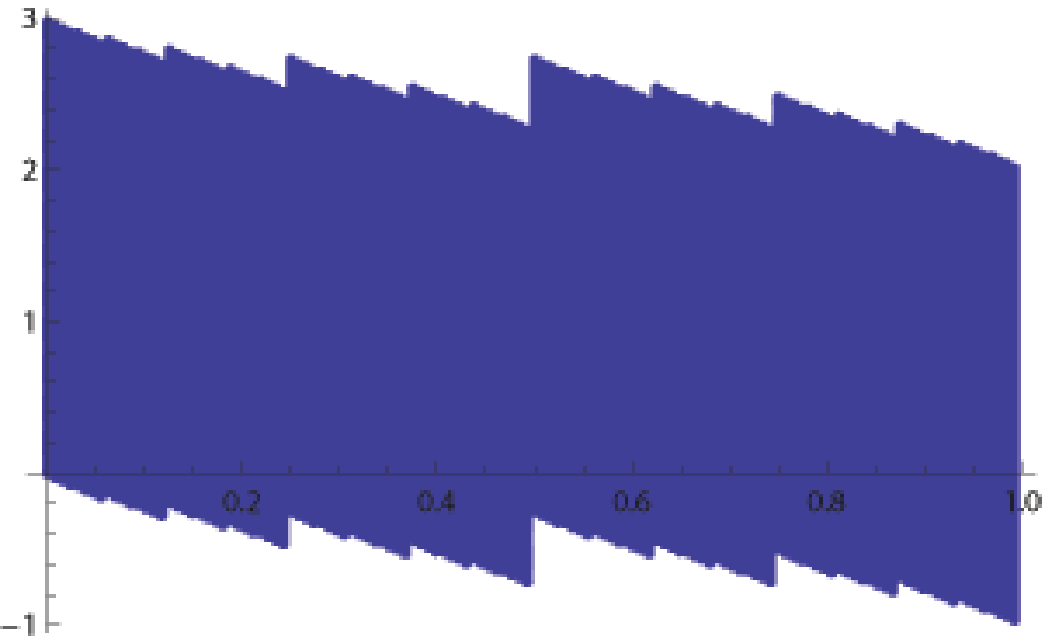}}} \hspace{3cm}
\vbox{\hbox{\epsfxsize 5cm\epsfbox{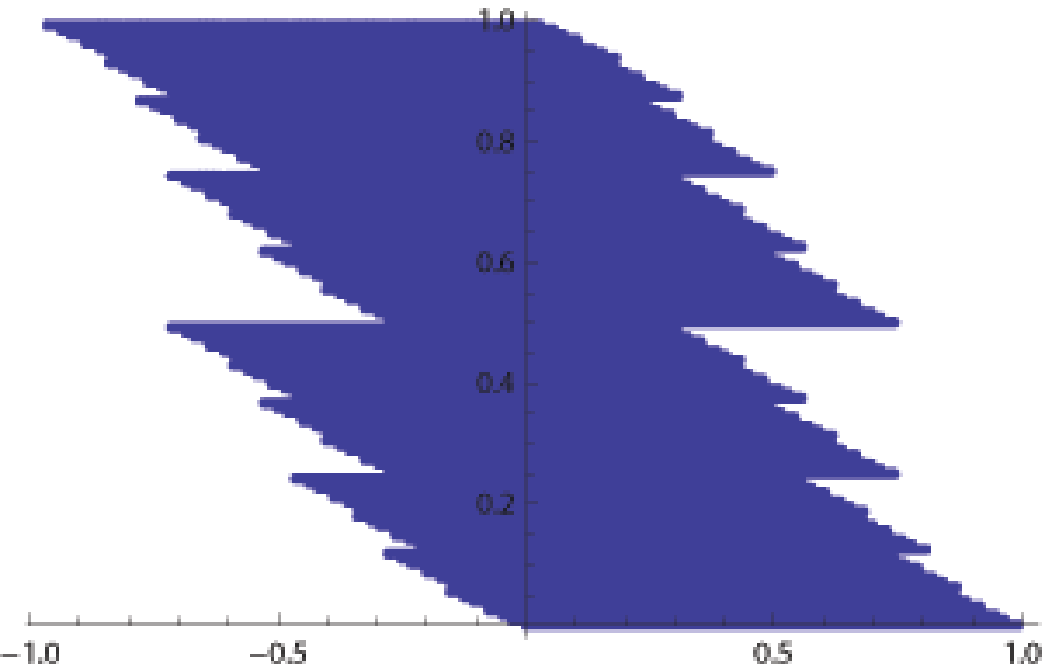}}}
}
\caption{ The attractors $X(B)$ and $X(L)$ for $p=3$}
\end{figure}

We will prove first that the invariant measure $\mu_B$ is the Lebesgue measure on the attractor $X_B$, renormalized so that the total measure is 1.

First, define the function $g:[0,1]\rightarrow\br$, for $x=\sum_{n=1}^\infty x_n/2^n$, with $x_n\in\{0,1\}$ for all $n$,
$$g(x)=-\sum_{n=1}^\infty\frac{n}{2}\frac{x_n}{2^n}.$$

The number $g(x)$ will represent the lowest endpoint of the vertical slice of $X_B$ at $x$. 
We claim that $g$ has the following property
\begin{equation}
g\left(\frac{x+b}{2}\right)=-\frac{x+b}{4}+\frac{g(x)}2,\quad(x\in[0,1],b\in \{0,1\})
\label{eqg1}
\end{equation}

Indeed, if $x=\sum_{n=1}^\infty x_n/2^n$ then $(x+b)/2=b/2+x_1/2^2+\dots+x_n/2^{n+1}+\dots$, so 
$$g\left(\frac{x+b}{2}\right)=-\left( -\frac12\frac b2+\frac22\frac{x_1}{2^2}+\dots+\frac{n+1}{2}\frac{x_n}{2^{n+1}}+\dots\right)$$
$$=-\left(\frac b4+\frac14\left(\frac{x_1}{2}+\dots+\frac{x_n}{2^n}+\dots\right)\right)+\frac12\left(\frac12\frac{x_1}{2}+\dots+\frac{n}{2}\frac{x_n}{2^n}+\dots\right)=-\frac{x+b}{4}+\frac{g(x)}{2}.$$

We prove now that the measure $\mu_B$ is given (up to the renormalization factor $1/p$) by 
\begin{equation}\label{eq4.7mu}
\int_0^1\int_0^p f(x,g(x)+y)\,dy\,dx.
\end{equation}
For this, we have to check the invariance equation. We have
$$\frac14\sum_{b_1\in\{0,1\}}\sum_{b_2\in\{0,p\}}\int_0^1\int_0^pf\circ\tau_{(b_1,b_2)^T}(x,g(x)+y)\,dy\,dx=$$$$
\frac14\sum_{b_1\in\{0,1\}}\sum_{b_2\in\{0,p\}}\int_0^1\int_0^p f\left(\frac{x+b_1}{2},-\frac{x+b_1}{4}+\frac{g(x)+y+b_2}{2}\right)\,dy\,dx$$
$$=\frac14\sum_{b_1\in\{0,1\}}\sum_{b_2\in\{0,p\}}\int_0^1\int_0^p f\left(\frac{x+b_1}{2},g\left(\frac{x+b_1}{2}\right)+\frac{y+b_2}2\right)\,dy\,dx$$
$$=\sum_{b_1\in\{0,1\}}\sum_{b_2\in\{0,p\}}\int_{\frac{b_1}2}^{\frac{b_1+1}{2}}\int_{\frac{b_2}{2}}^{\frac{b_2+p}{2}} f\left(x,g(x)+y\right)\,dy\,dx$$
$$=\int_0^1\int_0^p f\left(x,g(x)+y\right)\,dy\,dx.$$

\begin{lemma}\label{lem4.7.1}
The function $g$ is continuous at points $x$ which are dyadic-irrational (i.e., not of the form $x=k/2^n$), and it has jump discontinuities at dyadic rational values $x$. The sum of the jumps at the discontinuities is infinite. 
\end{lemma}

\begin{proof}
Set $x=1/2$, then since $1/2$ has two binary expansions $.0111\dots$ and $.1000\dots$, we get that 
$$g(1/2)=-1/4\mbox{ and }g(1/2)_-=-\sum_{n=2}^\infty\frac{n}{2^{n+1}}=-\frac34.$$
Therefore the jump at $x$ is $D(1/2)=1/2$. Using the relation \eqref{eqg1}, and induction we get that the jump at 
$k/2^n$ with $k$ odd, $k<2^n$ and $n\geq 1$ is $D(k/2^n)=1/2^{n+1}$.

Since there are $2^{n-1}$ distinct dyadics of the form $k/2^n$ with $0<k<2^n$, $k$ odd, the sum of the jumps is at least $\sum_{n=1}^\infty 2^{n-1}\frac1{2^n}=\infty$.

Since the binary representation is unique when $x$ is dyadic-irrational, it follows that $g$ is continuous at such points. 
\end{proof}

\begin{lemma}\label{lem4.7.3}
The function $g$ is nowhere differentiable. 

\end{lemma}

\begin{proof}
From the proof of Lemma \ref{lem4.7.1} we see that each interval of length $\frac{1}{2^n}$ contains a dyadic of the form $k/2^n$ so it contains a pair of points $x,y$ such that $|g(x)-g(y)|>\frac{1}{2^{n+2}}$.
Then 
$$\left|\frac{g(x)-g(y)}{x-y}\right|>\frac14.$$
\end{proof}

\begin{lemma}\label{lem4.7.4}
The two attractors $X(B)$ and $X(L)$  are given by (see Figure 1)
\begin{equation}
X(B)=\left\{(x,y)\in\br^2 : 0\leq x\leq 1, g(x)\leq y\leq g(x)+p\right\},
\label{eq4.7.7}
\end{equation}
\begin{equation}
X(L)=\left\{ (\xi,\eta)\in\br^2 : 0\leq \eta\leq 1, g(\eta)\leq \xi\leq g(\eta)+1\right\}.
\label{eq4.7.8}
\end{equation}
\end{lemma}

\begin{proof}
The result follows from \eqref{eq4.7mu} and the fact that $X(B)$ is the support of $\mu_B$. For $X(L)$ the proof is analogous. 
\end{proof}

The computation for \eqref{eq4.7mu} proves that 
\begin{lemma}\label{lem4.7.5}
The measure $\mu_B$ is the $\frac1p$ times the Lebesgue measure restricted to $X(B)$.
\end{lemma}

\begin{lemma}\label{lem4.7.6}
The Fourier transform of the measure $\mu_B$ is 
$$\hat\mu_B(t_1,t_2)=\int_0^1e^{2\pi i (t_1x+t_2g(x))}\,dx \cdot \frac{e^{2\pi i pt_2}-1}{2\pi ipt_2}$$
for $(t_1,t_2)\in\br^2$.
\end{lemma}

\begin{proof}
Immediate from \eqref{eq4.7mu}.
\end{proof}

\begin{lemma}\label{lem4.7.7}
Let $\mu_B$ as above. Then
\begin{equation}
\sum_{n\in\bz}|\hat\mu_B(t_1+n,t_2)|^2=\frac{\sin^2(p\pi t_2)}{(p\pi t_2)^2},\quad((t_1,t_2)\in\br^2)
\label{eq4.7.11}
\end{equation}
\end{lemma}

\begin{proof}

Using Lemma \ref{lem4.7.6} we get 
$$\sum_{n\in\bz}|\hat\mu_B(t_1+n,t_2)|^2=\sum_{n\in\bz}\left|\int_0^1e_n(x)e^{2\pi i(t_1x+t_2g(x))}\,dx\right|^2\frac{\sin^2(p\pi t_2)}{(p\pi t_2)^2}=$$
$$\quad( \mbox{ using Parseval's relation})\quad=\frac{\sin^2(p\pi t_2)}{(p\pi t_2)^2}.$$
\end{proof}

\begin{proposition}\label{prop4.7.8}
The set $\bz\times\frac1p\bz$ is a spectrum for the measure $\mu_B$. The set $X(B)$ tiles $\br^2$ by $\bz\times p\bz$. 
\end{proposition}

\begin{proof}
We have to prove that 
$$\sum_{\gamma\in\bz\times\frac1p\bz}|\hat\mu_B(t+\gamma)|^2=1,\quad(t\in\br^2).$$
(see e.g. \cite{DJ06b} for details). 

$$\sum_{\gamma\in\bz\times\frac1p\bz}|\hat\mu_B(t+\gamma)|^2=\sin^2(p\pi t_2)\sum_{k\in\bz}\frac{1}{\left( p\pi(t_2+\frac kp)\right)^2}$$ $$=\sin^2(p\pi t_2)\sum_{k\in\bz}\frac{1}{\pi^2(pt_2+k)^2}
=\frac{\sin^2(p\pi t_2)}{\sin^2(p\pi t_2)}=1.$$
\end{proof}

\begin{proposition}\label{pr4.7.8}
The $B$-extreme $L$-cycles are 
$$\left\{\begin{pmatrix}
0\\0
\end{pmatrix}
\right\}
,
\left\{\begin{pmatrix}
1\\0
\end{pmatrix}
\right\}
,
\left\{\begin{pmatrix}
0\\1
\end{pmatrix}
\right\}
,
\left\{\begin{pmatrix}
-1\\1
\end{pmatrix}
\right\}
$$
and the associated words are, respectively
$$\begin{pmatrix}
0\\0
\end{pmatrix},
\begin{pmatrix}
1\\0
\end{pmatrix},
\begin{pmatrix}
1\\1
\end{pmatrix},
\begin{pmatrix}
0\\1
\end{pmatrix}
$$
\end{proposition}

\begin{proof}
For this we can use Theorem \ref{th1.6}. Note that $\Gamma(B)^\circ=\bz\times\bz$.
\end{proof}

\begin{remark}\label{rem4.7.10}
For later use, we record the following points on the graph of $g$
\begin{equation}
g\left(\frac13\right)=-\frac49,\quad g\left(\frac23\right)=-\frac59.
\label{eq4.7.16}
\end{equation}
To see this, set $x=g(1/3)$ and $y=g(2/3)$. Use \eqref{eqg1}: clearly $\frac{\frac13+1}{2}=\frac23$ and $\frac{\frac23+0}{2}=\frac13$.
$$x=g\left(\frac{\frac23+0}{2}\right)=-\frac{\frac23+0}{4}+\frac y2$$
$$y=g\left(\frac{\frac13+1}{2}\right)=-\frac{\frac13+1}{4}+\frac x2$$
and then solve for $x$ and $y$.

The numbers in \eqref{eq4.7.16} result from intersecting the respective (horizontal) line segments $y = 1/3$ and $y = 2/3$ with the compact tile $X(L)$ in Figure 1. The invariant set $M$ (Definition \ref{def4.2}) is the union of these two segments; see also \eqref{eq4.7.8} in Lemma \ref{lem4.7.4}. Each line segment is unit-length, and the numbers in \eqref{eq4.7.16} are the respective left-hand side endpoints.

      The infinite random walk defined by restricting the maps  $\tau_l$  to $M$ ( Proposition \ref{pr4.2.1})  turns out in Example \ref{ex4.7} to be a zigzag motion, with the walker making successive jumps back and forth between the two line-segments.

      In the discussion below, we find the cycles $C$ from section 3, and the corresponding orthogonal sets  $\Lambda(C)$, (illustrated in Figure 2). Conclusion: Example \ref{ex4.7} has one non-cycle invariant set $M$, and 4 distinct cycles $C$.
\end{remark}

\begin{remark}\label{rem4.7.11}
We want to describe the sets $\Lambda(C)$ corresponding to the cycles in Proposition \ref{pr4.7.8}. We need some preliminaries. Each non-negative integer $k\geq0$ can be represented in base 2 as 
$$k=j_0+2j_1+2^2j_2+\dots+ 2^nj_n.$$
The representation is unique, up to adding $0$ digits at the end. 
Each negative integer $k<0$ can be represented as
$$k=-2^{n+1}+j_0+2j_1+2^2j_2+\dots+2^nj_n.$$
The representation is unique up to adding $1$ digits at the end. 

We define the following function:
\begin{equation}
\mbox{ For }k\geq 0, k=j_0+2j_1+\dots+2^nj_n,\quad h(k):=j_1+2\cdot 2^1j_2+\dots+n2^{n-1}j_n,
\label{eq4.7h1}
\end{equation}
\begin{equation}
\mbox{ For }k<0, k=-2^{n+1}+j_0+2j_1+\dots+2^nj_n,\quad h(k):=2^{n+1}-(n+1)2^n+j_1+2\cdot 2^1j_2+\dots+n2^{n-1}j_n.
\label{eq4.7h2}
\end{equation}
It is easy to check that the function $h$ is well defined, i.e., it does not depend on the representation, and it satisfies the following equation
\begin{equation}
2h(n)+n=h(2n+j_0),\mbox{ for all }n\in\bz, j_0\in\{0,1\}.
\label{eq4.7h3}
\end{equation}
\end{remark}

\begin{proposition}\label{pr4.7.11}
The four sets $\Lambda(C)$ (as in Theorem \ref{th1.8}) for the four cycles in Proposition \ref{pr4.7.8} are 
$$\Lambda\left(\left\{\begin{pmatrix}0\\0\end{pmatrix}\right\}\right)=\left\{\begin{pmatrix}t_1\\t_2\end{pmatrix}\in\bz^2 : t_2\geq 0, t_1\geq h(t_2)\right\},$$
$$\Lambda\left(\left\{\begin{pmatrix}-1\\1\end{pmatrix}\right\}\right)=\left\{\begin{pmatrix}t_1\\t_2\end{pmatrix}\in\bz^2 : t_2< 0 , t_1\geq h(t_2)\right\},$$
$$\Lambda\left(\left\{\begin{pmatrix}0\\1\end{pmatrix}\right\}\right)=\left\{\begin{pmatrix}t_1\\t_2\end{pmatrix}\in\bz^2 : t_2< 0, t_1< h(t_2)\right\},$$
$$\Lambda\left(\left\{\begin{pmatrix}1\\0\end{pmatrix}\right\}\right)=\left\{\begin{pmatrix}t_1\\t_2\end{pmatrix}\in\bz^2 : t_2\geq 0, t_1<h(t_2)\right\},$$
and as a result 
$$\bigcup_{C}\Lambda(C)=\bz^2,\quad \Lambda(0,0)\cup\Lambda(1,0)=\bz\times\bz_{\geq0}.$$
\end{proposition}

\begin{figure}[ht]\label{fig2}
\centerline{ 
\vbox{\hbox{\epsfxsize 5cm\epsfbox{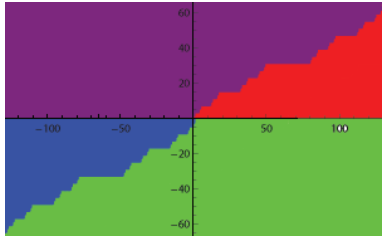}}}}
\caption{ The sets $\Lambda(C)$}
\end{figure}

\begin{proof}
Let $A=R^T$. Take a point $(k,l)\in\bz^2$. Then, since $L$ is a complete set of representatives for $\bz^2/A\bz^2$, it follows that there is a unique $(m,n)\in\bz^2$ and $(l_0,l_1)\in L$ such that 
$(k,l)^T=A(m,n)^T+(l_0,l_1)^T$. We consider the map $r:(k,l)\mapsto (m,n)$ Using \eqref{eq4.7h3}, a computation shows that the four sets $\Lambda(C)$ are invariant under this map. Iterating the map $r$, each point is eventually maped into $-C$ for one of the cycles $C$. For more details check also \cite{BrJo99, DuJo09e}.   
\end{proof}
\begin{remark}
The sum of the subspaces $\H(C)$ (as in Theorem \ref{th1.8}) where $C$ runs over the $B$-extreme $L$-cycles is not total in $L^2(\mu_B)$. This is because the spectrum of $\mu_B$ is $\bz\times \frac1p\bz$, not $\bz\times\bz$. Actually, all the exponential functions $e_{(t,l/p)^T}$, $l\not\equiv0\mod p$  are orthogonal to this sum. 
Indeed, for $(n_1,n_2)\in\bz^2$:
$$\ip{e_{(t,l/p)^T}}{e_{(n_1,n_2)^T}}_{L^2(\mu_B)}=\int_0^1e^{2\pi i((t-n_1)x+(\frac lp-n_2)g(x))}\int_0^pe^{2\pi i(\frac lp-n_2)y}\,dy\,dx=0.$$
\end{remark}

\begin{lemma}\label{lem4.7.12}
Let $(R,B,L)$ be a Hadamard system in $\br^d$, and let $\Gamma\subset\br^d$ be a subset such that the vectors $\{e_\gamma: \gamma\in\Gamma\}$ are orthogonal $L^2(\mu_B)$. Let $P_\Gamma$ be the projection onto the corresponding subspace in $L^2(\mu_B)$ spanned by these vectors, and let 
$$h_{P_\Gamma}(t):=\|P_{\Gamma}e_{-t}\|^2,\quad(t\in\br^d).$$
Then 
$$h_{P_\Gamma}(t)=\sum_{\gamma\in\Gamma}|\hat\mu_B(t+\gamma)|^2,\quad(t\in\br^d)$$
and $0\leq h_{P_\Gamma}(t)\leq 1$ for all $t\in\br^d$.
\end{lemma}

\begin{proof}
Since the vectors are ortogonal we have 
\begin{equation}\label{eq4.7.20}
P_\Gamma=\sum_{\gamma\in\Gamma}\left| e_\gamma\rangle\langle e_\gamma\right|
\end{equation}
with use of Dirac notation for rank-1 operators. Since 
$$\left\|\left|e_\gamma\rangle\langle e_\gamma\right| e_{-t}\right\|^2=|\hat\mu_B(t+\gamma)|^2,$$
a substitution into \eqref{eq4.7.20} yields the desired conclusions. 
\end{proof}

\begin{corollary}\label{cor4.7.13}
For Example \ref{ex4.7}, corresponding to the four cycles in Proposition \ref{pr4.7.11}, we get the four independent harmonic functions 
\begin{equation}
h_C(t)=\left\|P_{\Lambda(C)}e_{-t}\right\|^2,\quad(t\in\br^2),
\label{eq4.7.21}
\end{equation}
and
\begin{equation}
\sum_{C}h_C(t)=\left\|P_{\bz^2}e_{-t}\right\|^2=\left(\frac{\sin p\pi t_2}{p\sin\pi t_2}\right)^2,\quad((t_1,t_2)\in\br^2)
\label{eq4.7.22}
\end{equation}

\end{corollary}

\begin{proof}
It follows from Proposition \ref{pr4.7.11} that the sum of the four harmonic functions in \eqref{eq4.7.22} is $\|P_{\bz^2}e_{-t}\|^2$. An application of Lemma \ref{lem4.7.12} and \ref{lem4.7.7} then yields
$$\|P_{\bz^2}e_{-t}\|^2=\frac{\sin^2(p\pi t_2)}{(p\pi)^2}\sum_{n\in\bz}\frac{1}{(t+n)^2}=\frac{\sin^2(p\pi t_2)}{p^2\sin^2(\pi t_2)}$$
which is the desired conclusion \eqref{eq4.7.22}.
\end{proof}

\begin{remark}\label{rem4.7.19}
The functions 
$$f_A(t)=\ip{Ae_{-t}}{e_{-t}}_{L^2(\mu_B)}$$ 
in Lemma \ref{lem4.7.12} are entire analytic on $\br^d$ for every bounded operator $A$. In fact the functions in Corrolary \ref{cor4.7.13} are trigonometric polynomials. Indeed, the expression on the right hand side in \eqref{eq4.7.22} divides, for example, for $p=3$, we have 
$$\frac{\sin 3\pi t}{\sin \pi t}=3-4\sin^2(\pi t)=4\cos^2(\pi t)-1.$$
\end{remark}

\begin{remark}\label{rem4.7.20}
In affine examples, it is typically more difficult to obtain explicit formulas for the individual harmonic functions $h_C(t)$ (as in Proposition \ref{pr4.7.11}) than for their summations. The difference may be explained by the difference between number theoretic expressions such as 
\begin{equation}
\sum_{n\in\bz}\frac{1}{(t+n)^2}=\frac{\pi^2}{\sin^2(\pi t)}
\label{eq20.1}
\end{equation}
on the one hand, and 
\begin{equation}
\zeta_2(t)=\sum_{n=0}^\infty\frac{1}{(t+n)^2}=\left(\frac{d}{dt}\right)^2\log \Gamma(t),
\label{eq20.2}
\end{equation}
i.e., summation over $\bz_{\geq0}$, where $\Gamma$ is the usual Gamma function. The function $\zeta_2$ in \eqref{eq20.2} is one of the Hurwitz-Riemann zeta functions \cite{AAR99}. 

Using this, we arrive at the following:
$$h_{\Lambda(0,0)}(t)+h_{\Lambda(1,0)}(t)=\frac{\sin^2(p\pi t_2)}{\pi^2}\zeta_2(t_2)$$
and
$$h_{\Lambda(-1,1)}(t)+h_{\Lambda(0,1)}(t)=\frac{\sin^2(p\pi t_2)}{\pi^2}\left(-\frac{1}{t_2^2}+\zeta_2(-t_2)\right)$$
for all $t=(t_1,t_2)\in\br^2$.
In all cases the functions $t\mapsto h_{\Lambda(C)}(t)$ have entire analytic extensions.
\end{remark}

\begin{proposition}\label{pr4.7.15}
Fi $p\in\bz$ odd, and let $(R,B,L)$ be as in Example \ref{ex4.7}. List the non-trivial orbits $O_1,O_2,\dots$ for $x\mapsto 2x$ acting on $\bz_p=\bz/p\bz=\{0,1,2,\dots,p-1\}$. Then the sets 
\begin{equation}
M(O_i)=X(L)\cap\left\{\begin{pmatrix} \xi\\ \eta\end{pmatrix} : p\eta\in O_i\right\}
\label{eq4.7.15}
\end{equation}
are minimal invariant sets 
\end{proposition}

\begin{remark}\label{rem4.7.16}
If $p=3$, there is only one orbit $O_1=\{1,2\}$. For $p=5$: $O_1=\{1,2,4,3\}$. If $p=7$, the list of non-trivial orbits is $O_1=\{1,2,4\}$, and $O_2=\{3,6,5\}$. If $p=9$ the list is 
$O_1=\{1,2,4,8,7,5\}$, $O_2=\{3,6\}$. 

For general values of $p$, the orbits can be computed, and the number of possible orbits gets arbitrarily large.
\end{remark}

\begin{proof}[Proof of Proposition \ref{pr4.7.15}]
This is based on the formula \eqref{eq4.1} in the proof of Theorem \ref{pr4.4} above. When applied to Example \ref{ex4.7}, we get 
\begin{equation}\label{eq4.7.26}
\tau_{(l_1,l_2)}^{(L)}(t_1,t_2)=\begin{pmatrix}
\frac12(t_1+l_1-\frac{t_2+l_2}{2})\\ \frac{t_2+l_2}{2}\end{pmatrix}
\end{equation}
Then 
\begin{equation}
\left|\chi_B\left(\tau_{(l_1,l_2)}^{(L)}(t_1,t_2)\right)\right|^2=\cos^2(\pi\frac12(t_1+l_1-\frac{t_2+l_2}2)\cos^2(\frac{p\pi}{2}(t_2+l_2)),
\label{eq4.7.27}
\end{equation}
and the conclusion follows.

The details for $p=3$ are as follows ($l\in\{0,1\}$):

$$\left|\chi_B\left(\tau_{(l,1)}^{(L)}(t,1/3)\right)\right|^2=\cos^2(\frac{\pi}2(t+l-\frac23))$$
while 
$$\left|\chi_B\left(\tau_{(l,0)}^{(L)}(t,1/3)\right)\right|^2=0.$$

The transition probabilities in the other direction are 
$$\left|\chi_B\left(\tau_{(l,0)}^{(L)}(t,2/3)\right)\right|^2=\cos^2(\frac{\pi}2(t+l-\frac13)),$$
and
$$\left|\chi_B\left(\tau_{(l,1)}^{(L)}(t,2/3)\right)\right|^2=0.$$
With the use of \eqref{eq4.7.26} and \eqref{eq4.7.27} we get the remaining conclusion of the proposition.
\end{proof}

\begin{corollary}\label{cor4.7.17}
Let $p=3$ in Example \ref{ex4.7}. Then there is just one non-cycle minimal invariant set 
$$M=X(L)\cap\left\{\begin{pmatrix}\xi \\ \eta\end{pmatrix} : \eta\in\left\{\frac13,\frac23\right\}\right\}=\left\{\begin{pmatrix}t \\ \frac13\end{pmatrix} : -\frac49\leq t\leq\frac 59\right\}\cup
\left\{\begin{pmatrix}t \\ \frac23\end{pmatrix} : -\frac59\leq t\leq\frac 49\right\}.$$

If $Q$ is the projection onto $\H(M)\subset L^2(\mu_B)$, then the associated harmonic function is :
\begin{equation}
h_Q(t)=1-(1-\frac43\sin^2(\pi t_2))^2,\quad (t=(t_1,t_2)\in\br^2)
\label{eq4.7.28}
\end{equation}
\end{corollary}

\begin{proof}
Combine the results in Lemma \ref{lem4.7.12}, Proposition \ref{pr4.7.15} and Remark \ref{rem4.7.10}.

\end{proof}

\end{example}

\bibliographystyle{alpha}
\bibliography{cuntz}

\end{document}